\documentclass[12pt]{article}

\usepackage[english,russian]{babel}
\usepackage{latexsym,amssymb,mathrsfs,amsfonts}

\newcommand{\Lam}{{\cal L}}

\newcommand{\Ss}{{\cal S}}
\newcommand{\N}{{\cal N}}
\newcommand{\B}{{\cal B}}

\begin{document}

{\bf Sets of vector fields with various shadowing properties}
\bigskip

S. Yu. Pilyugin and S. B. Tikhomirov
\footnote{Saint Petersburg State University,
199164, St.Petersburg, Universitetskaya nab., 7/9.}

\bigskip

Let $M$ be a smooth closed manifold with 
Riemannian metric $\mbox{dist}$.
We consider the space of tangent vector fields on
$M$ of class $C^1$ with the  $C^1$ topology;
denote by $\phi(t,x)$ the trajectory of a field
$X$ such that $\phi(0,x)=x$.

Fix a number $d>0$. We say that a mapping $g:\mathbb{R}\to M$
is a $d$-{\em pseudotrajectory} of a flow $\phi$
(and of the corresponding field $X$) if the inequalities
$$
\mbox{dist}(\phi(t,g(\tau)),g(\tau+t))<d,\quad |t|\leq 1,
$$
hold for any $\tau\in\mathbb{R}$.

The shadowing problem is related to the following question:
under which condition, for any pseudotrajectory of a field $X$
there exists a close trajectory? The study of this problem 
was originated by D. V. Anosov [1] and R. Bowen [2];
the modern state of the shadowing theory is reflected in the
monographs [3, 4].

Let us note that the main difference between the shadowing
problem for flows and the similar problem for discrete dynamical
systems generated by diffeomorphisms is related to the necessity
of reparametrization of shadowing trajectories in the former
case.

The aim of this short note is to describe the structure of
$C^1$-interiors of sets of vector fields with various shadowing
properties.

A monotonically increasing homeomorphism $h$ of the line 
$\mathbb{R}$ such that $h(0)=0$ is called a
{\em reparametrization}. 

Let $a>0$; denote by $\mbox{Rep}(a)$ the set of 
reparametrizations $h$ such that
$$
\left|\frac{h(t_1)-h(t_2)}{t_1-t_2}-1\right|<a
$$
for any different $t_1,t_2\in\mathbb{R}$.

Let us define the main shadowing properties which we study
in this paper.

We say that a field $X$ has the  {\em regular shadowing property}
if for any $\epsilon>0$ there exists a number $d>0$ with the
following property: for any $d$-pseudotrajectory $g$ of the field $X$ 
there exists a point $p\in M$ and a reparametrization $h\in\mbox{Rep}(\epsilon)$
such that
\begin{equation}
\label{1}
\mbox{dist}(\phi(h(t),p),g(t))<\epsilon,\quad t\in\mathbb{R}.
\end{equation}
Denote by RegSh the set of vector fields having the 
regular shadowing property.

We say that a field $X$ has the  {\em Lipschitz shadowing property}
if there exist numbers $d_0,\Lam>0$ having the following property:
for any $d$-pseudotrajectory $g$ of the field $X$ with $d\leq d_0$ 
there exists a point $p\in M$ and a reparametrization $h\in\mbox{Rep}(\Lam d)$
such that inequalities (\ref{1}) hold with $\epsilon=\Lam d$.
Denote by LipSh the set of vector fields having the 
Lipschitz shadowing property.

We say that a field $X$ has the {\em oriented shadowing property}
if for any $\epsilon>0$ there exists a number $d>0$ having the following property:
for any $d$-pseudotrajectory $g$ of the field $X$ there exists
a point $p\in M$ and a reparametrization $h$ such that inequalities (\ref{1})
hold (thus, we do not require the reparametrization $h$ 
to be close to the identity).
Denote by OrientSh the set of vector fields having the 
Lipschitz shadowing property.

Finally, we say that a field $X$ has the {\em orbital shadowing property}
if for any  $\epsilon>0$ there exists a number $d>0$ having the following property:
for any $d$-pseudotrajectory $g$ of the field $X$ there exists
a point $p\in M$ such that
$$
\mbox{dist}_H(\mbox{Cl}\{\phi(t,p):t\in\mathbb{R}\},
\mbox{Cl}\{g(t):t\in\mathbb{R}\})<\epsilon,
$$
where $\mbox{Cl}A$ is the closure of a set $A$, and $\mbox{dist}_H$ 
is the Hausdorff distance.
Denote by OrbitSh the set of vector fields having the 
orbital shadowing property.

Clearly, the following inclusions hold:
$$
\mbox{LipSh}\subset \mbox{RegSh}\subset \mbox{OrientSh}\subset
\mbox{OrbitSh}.
$$

We introduce the following notation: $\Ss$ denotes the set of
structurally stable vector fields, and $\N$ denotes the set of
nonsingular vector fields. If $P$ is a set of vector fields, we
denote by $\mbox{Int}^1(P)$ the interior of the set $P$
with respect to the $C^1$ topology.

It was shown in [5] that $\Ss\subset\mbox{LipSh}$.

We define the following class of vector fields which is important
for us. We say the a field $X$ belongs to the class $\B$
if it has hyperbolic rest points $p$ and $q$ (not necessarily distinct)
such that

(1) the Jacobi matrix $DX(p)$ has a pair of complex conjugate
eigenvalues $a_1\pm b_1 i$ with $a_1<0$ of multiplicity 1, and if $c_1+d_1 i$ 
is an eigenvalue different from $a_1\pm b_1 i$ and such that $c_1<0$, then $c_1<a_1$;

(2) the Jacobi matrix $DX(q)$ has a pair of complex conjugate
eigenvalues $a_2\pm b_2 i$ with $a_2>0$ of multiplicity 1, and if $c_2+d_2 i$ 
is an eigenvalue different from $a_2\pm b_2 i$ and such that $c_2>0$, then $c_2>a_2$;

(3) the unstable manifold $W^u(p)$ and the stable manifold $W^s(q)$
have a trajectory of nontransverse intersection.
\medskip

{\bf Theorem 1. } $\mbox{Int}^1(\mbox{OrbitSh})\cap\N\subset\Ss.$
\medskip

This theorem generalizes the main result of the recent paper [6]
where it was shown that $\mbox{Int}^1(\mbox{RegSh})\cap\N\subset\Ss$. 
For discrete dynamical systems generated by diffeomorphisms,
an analog of Theorem 1 was obtained in [7].
\medskip

{\bf Theorem 2. } $\mbox{Int}^1(\mbox{OrientSh}\setminus\B)\subset\Ss.$
\medskip

{\bf Theorem 3. } $\mbox{Int}^1(\mbox{LipSh})\subset\Ss.$
\medskip

The above-mentioned result of [5] and Theorem 3 imply that
$$
\mbox{Int}^1(\mbox{LipSh})=\Ss.
$$

{\bf Theorem 4. } {\em If }$\mbox{dim}M\leq 3$, {\em then } $\mbox{Int}^1(\mbox{OrientSh})\subset\Ss.$
\medskip

Thus, if $\mbox{dim}M\leq 3$, then
$$
\mbox{Int}^1(\mbox{OrientSh})=\Ss.
$$
\medskip

References
\medskip

1. Anosov, D.V. // {\em Trudy 5th Int. Conf. Nonlin. Oscill.}, Kiev: 1970, vol. 2, pp. 39-45.

2. Bowen, R. // {\em Lect. Notes Math.}, 1975, vol. 470.

3. Pilyugin, S.Yu. // {\em Lect. Notes Math.}, 1999,  vol. 1706.

4. Palmer, K. {\em Shadowing in Dynamical Systems. Theory and Applications},
Dordrecht-Boston-London: Kluwer Acad. Publ., 2000.

5. Pilyugin, S.Yu. // {\em J. Differ. Equations}, 1997, vol. 140, pp. 238-265.

6. Lee, K. and Sakai, K. // {\em J. Differ. Equations}, 2007. vol. 232. pp. 303-313.

7. Pilyugin, S.Yu., Rodionova, A.A., and Sakai K. // {\em Discrete Contin. 
Dyn. Syst.}, 2003. vol. 9, pp. 287-308.
\medskip

Translated by S. Yu. Pilyugin.

\end{document}